\documentclass[graybox]{svmult}

\usepackage{helvet}         
\usepackage{courier}        
\usepackage{type1cm}        
\usepackage{makeidx}         
\usepackage{graphicx}        
\usepackage{multicol}        
\usepackage[bottom]{footmisc}
\usepackage{eufrak}
\usepackage{color}
\begin{document}

\title*{Quadratic relations of\\
the deformed $W$-algebra}
\author{Takeo Kojima}
\institute{Takeo Kojima\at
Yamagata University, Jonan 4-chome 3-16, Yonezawa 992-8510, JAPAN\\
\email{kojima@yz.yamagata-u.ac.jp}}
%
%
\maketitle

\abstract{
The deformed $W$-algebra is a quantum deformation of 
the $W$-algebra ${\cal W}_\beta(\mathfrak{g})$ in conformal field theory.
Using the free field construction, we obtain a closed set of quadratic relations of the $W$-currents
of the deformed $W$-algebra.
This allows us to define the deformed $W$-algebra by generators and relations.
In this review,
we study two types of deformed $W$-algebra. 
One is the deformed $W$-algebra 
${\cal W}_{x,r}\big(A_{2N}^{(2)}\big)$,
and the other is the $q$-deformed corner vertex algebra 
$q$-$Y_{L_1, L_2, L_3}$ that is a 
generalization of the deformed $W$-algebra
${\cal W}_{x,r}\big(A(M,N)^{(1)}\big)$ 
via the quantum toroidal algebra.}

\section{Introduction}
\label{sec:1}
The deformed $W$-algebra 
${\cal W}_{x,r}({\mathfrak g})$ is both a two-parameter deformation of the classical $W$-algebra
${\cal W}({\mathfrak g})$ in soliton theory
and a one-parameter deformation of the $W$-algebra 
${\cal W}_\beta({\mathfrak g})$ 
in conformal field theory.
The deformation theory of 
the $W$-algebra ${\cal W}_\beta({\mathfrak g})$ 
has been studied in papers \cite{
Shiraishi-Kubo-Awata-Odake,
Awata-Kubo-Odake-Shiraishi,
Feigin-Frenkel,
Brazhnikov-Lukyanov,
Frenkel-Reshetikhin1,
Sevostyanov,
Odake,
Feigin-Jimbo-Mukhin-Vilkoviskiy,
Ding-Feigin,
Kojima1,
Kojima2,
Kojima3,
Harada-Matsuo-Noshita-Watanabe}.
In comparison with the conformal case, the deformation theory of the $W$-algebra is still not fully understood.
Except in low-rank cases such as 
the Virasoro algebra and the ${\cal W}_3$-algebra, 
it isn't easy to handle 
the $W$-algebras ${\cal W}_\beta({\mathfrak g})$ 
in a computational way \cite{Kodera-Ueda}.
In the case of the deformed $W$-algebra, 
it is sometimes 
possible to perform concrete calculations relatively easily.
For instance, quadratic relations of the deformed $W$-algebra 
${\cal W}_{x,r}(\mathfrak{g})$
have already been known in the cases of 
$\mathfrak{g}=A_N^{(1)}$ and $A_{2}^{(2)}$
\cite{
Shiraishi-Kubo-Awata-Odake,
Awata-Kubo-Odake-Shiraishi,
Feigin-Frenkel,
Brazhnikov-Lukyanov, Odake}.
In the case of ${\cal W}_{x,r}\big(A_{1}^{(1)}\big)$,
the basic $W$-current ${\mathbf T}_1(z)$ satisfy
the following quadratic relation \cite{Shiraishi-Kubo-Awata-Odake}
\begin{eqnarray}
&&
g\left(\frac{z_2}{z_1}\right)
{\mathbf T}_1(z_1){\mathbf T}_1(z_2)-
g\left(\frac{z_1}{z_2}\right)
{\mathbf T}_1(z_2){\mathbf T}_1(z_1)
=c \left(
\delta\left(\frac{x^{-2}z_2}{z_1}\right)
-\delta\left(\frac{x^{2}z_2}{z_1}\right)
\right)
\nonumber
\end{eqnarray}
with an appropriate constant 
$c$ and a function $g(z)$.
In the case of ${\cal W}_{x,r}\big(A_{2}^{(2)}\big)$,
the basic $W$-current ${\mathbf T}_1(z)$ satisfy
the following quadratic relation \cite{Brazhnikov-Lukyanov}
\begin{eqnarray}
&&
f\left(\frac{z_2}{z_1}\right)
{\mathbf T}_1(z_1){\mathbf T}_1(z_2)-
f\left(\frac{z_1}{z_2}\right)
{\mathbf T}_1(z_2){\mathbf T}_1(z_1)\nonumber\\
&=&\delta\left(\frac{x^{-2}z_2}{z_1}\right)
{\mathbf T}_1(x^{-1}z_2)
-\delta\left(\frac{x^{2}z_2}{z_1}\right)
{\mathbf T}_1(xz_2)+c \left(
\delta\left(\frac{x^{-3}z_2}{z_1}\right)
-\delta\left(\frac{x^{3}z_2}{z_1}\right)
\right)
\nonumber
\end{eqnarray}
with an appropriate constant 
$c$ and a function $f(z)$.
In this review, the author 
would like to report 
the quadratic relations 
in the cases of the twisted algebra 
${\cal W}_{x,r}\big(A_{2N}^{(2)}\big)$ \cite{Kojima3} and 
the $q$-deformed corner vertex algebra 
$q$-$Y_{L_1, L_2, L_3}$ 
that is a generalization of the deformed $W$-algebra
${\cal W}_{x, r}\big(A(M,N)^{(1)}\big)$
via the quantum toroidal algebra
\cite{Feigin-Jimbo-Mukhin-Vilkoviskiy, Kojima1, Kojima2, 
Harada-Matsuo-Noshita-Watanabe}.
These relations allow us to define 
the deformed $W$-algebras by generators and relations.

The text is organized as follows.
In Section 2,
we review the quantum toroidal algebra ${\cal E}$ 
associated to $\mathfrak{gl}_1$
and the quantum algebra ${\cal K}$.
In Section 3,
we review the free field constructions
of the basic $W$-currents ${\mathbf T}_1(z)$
both for ${\cal W}_{x,r}\big(A_{2N}^{(2)}\big)$ and 
$q$-$Y_{L_1, L_2, L_3}$.
We introduce 
the higher $W$-currents ${\mathbf T}_i(z)$, $i=2,3,4,\ldots$, 
by
fusion procedure.
We present a closed set of quadratic relations.
Using these relations, 
we define the deformed $W$-algebras by generators and relations.

\section{Quantum toroidal algebra ${\cal E}$ associated to $\mathfrak{gl}_1$}
\label{sec:2}

\subsection{Notation}
Throughout the text we fix three complex parameters 
$q_1, q_2, q_3 \in {\mathbf C}^{\times}$
such that $q_1q_2q_3=1$.
We assume $q_1^l q_2^m q_3^n=1$ $(l, m, n \in {\mathbf Z})$ implies $l=m=n=0$.
We use the notation
$
s_c=q_c^{\frac{1}{2}}~(c=1,2,3), 
~\kappa_r=\prod_{c=1}^3(1-q_c^r)~(r \in {\mathbf Z}).
$
For any integer $n$, define $q$-integer 
\begin{eqnarray}
~[n]_q=\frac{q^n-q^{-n}}{q-q^{-1}}\nonumber
\end{eqnarray}
for complex number $q\neq 0$.
We use symbols or infinite products
\begin{eqnarray}
(a; p)_\infty=
\prod_{k=0}^\infty (1-a p^k),~~~
(a_1, a_2, \ldots, a_N; p)_\infty=
\prod_{i=1}^N (a_i;p)_\infty
\nonumber
\end{eqnarray}
 for $|p|<1$ and $a, a_1,\ldots, a_N \in {\mathbf C}$.
Define $\delta(z)$ by the formal series
$$\delta(z)=\sum_{m \in {\mathbf Z}}z^m.$$

\subsection{Quantum toroidal algebra ${\cal E}$
associated to ${\mathfrak{gl}}_1$}

In this section, we review the quantum toroidal algebra ${\cal E}$ associated to 
${\mathfrak{gl}}_1$
in Refs.
\cite{Feigin-Jimbo-Mukhin-Vilkoviskiy, Ding-Iohara, Feigin-Jimbo-Miwa-Mukhin}.
We set
\begin{eqnarray}
&&
g(z,w)=\prod_{j=1}^3(z-q_jw),~~~
\bar{g}(z,w)=\prod_{j=1}^3(z-q_j^{-1}w).
\nonumber
\end{eqnarray}
The quantum toroidal algebra ${\cal E}$ associated to $\mathfrak{gl}_1$
is an associative algebra with parameters $q_1, q_2, q_3$
generated by $e_n, f_n$ $(n \in {\mathbf Z})$, $h_r$ $(r \in {\mathbf Z}_{r \neq 0})$
and invertible central element $C$.
We set the currents
\begin{eqnarray}
e(z)=\sum_{n \in {\mathbf Z}}e_n z^{-n},
~f(z)=\sum_{n \in {\mathbf Z}}f_nz^{-n},
~\psi^\pm(z)=\exp\left(\sum_{r>0}\kappa_r h_{\pm r}z^{\mp r}\right).
\nonumber
\end{eqnarray}
The defining relations are given as follows.
\begin{eqnarray}
&&~[h_r, h_s]=\delta_{r+s,0}\frac{1}{\kappa_r}
\frac{C^r-C^{-r}}{r},
\nonumber
\\
&&
g(z,w)
\psi^+(C^{-1}z)e(w)=
\bar{g}(z,w)e(w)\psi^+(C^{-1}z),
\nonumber\\
&&
g(z,w)
\psi^-(z)e(w)=
\bar{g}(z,w)e(w)\psi^-(z),
\nonumber\\
&&
\bar{g}(z,w)
\psi^+(z)f(w)={g}(z,w)f(w)\psi^+(z),
\nonumber\\
&&
\bar{g}(z,w)
\psi^-(C^{-1}z)f(w)=
{g}(z,w)f(w)\psi^-(C^{-1}z),
\nonumber\\
&&
~[e(z), f(w)]=
\frac{1}{\kappa_1}\left(\delta\left(\frac{Cw}{z}\right)\psi^+(w)-
\delta\left(\frac{Cz}{w}\right)\psi^-(z)\right),\nonumber\\
&&g(z,w)e(z)e(w)=\bar{g}(z,w)e(w)e(z),~~
\bar{g}(z,w)f(z)f(w)={g}(z,w)f(w)f(z),
\nonumber
\\
&&{\rm Sym}_{z_1, z_2, z_3}\frac{z_2}{z_3}[e(z_1), [e(z_2), e(z_3)]]=0,
~~{\rm Sym}_{z_1, z_2, z_3}\frac{z_2}{z_3}[f(z_1), [f(z_2), f(z_3)]]=0, \nonumber
\end{eqnarray}
where we used
${\rm Sym}_{z_1, z_2, z_3}F(z_1, z_2, z_3)=
\frac{1}{3!}\sum_{\sigma \in S_3}
F(z_{\sigma(1)}, z_{\sigma(2)}, z_{\sigma(3)})$.
The quantum toroidal algebra ${\cal E}$ is endowed with a topological
Hopf algebra structure 
$({\cal E}, \Delta, \varepsilon, S)$.
We define the topological coproduct 
$\Delta: {\cal E}\to {\cal E} \tilde{\otimes} {\cal E}$, 
the counit $\varepsilon: {\cal E} \to {\mathbf C}$, 
and the antipode $S: {\cal E} \to {\cal E}$ 
as follows.
\begin{eqnarray}
&&
\Delta e(z)=e(C_2^{-1}z)\otimes \psi^+(C_2^{-1}z)+1\otimes e(z),\nonumber
\\
&&
\Delta f(z)=f(z)\otimes 1+\psi^-(C_1^{-1}z)\otimes f(C_1^{-1}z),\nonumber
\\
&&\Delta \psi^+(z)=\psi^+(z)\otimes \psi^+(C_1z),
\nonumber
\\
&&\Delta \psi^-(z)=\psi^-(C_2 z)\otimes \psi^-(z),~~
\Delta (C)=C \otimes C,\nonumber
\end{eqnarray}
where $C_1=C\otimes 1$, $C_2=1\otimes C$.
\begin{eqnarray}
&&
\varepsilon(e(z))=0,~~
\varepsilon(f(z))=0,~~
\varepsilon(\psi^\pm(z))=1,~~
\varepsilon(C)=1,\nonumber\\
&&
\tilde{e}(z)
=S(e(z))=-e(Cz)\psi^+(z)^{-1},\nonumber
\\
&&
\tilde{f}(z)
=S(f(z))=-\psi^-(z)^{-1}f(Cz),\nonumber
\\
&&
S(\psi^\pm(z))=\psi^\pm(C^{-1}z),~~
S(C)=C^{-1}.\nonumber
\end{eqnarray}

The quantum toroidal algebra ${\cal E}$ has three families of Fock representations
${\cal F}_c(u)$,
where $c=1,2,3$ and $u \in {\mathbf C}^{\times}$.
We call $c$ the color.
The Fock module ${\cal F}_c(u)$ has level $s_c$.
The Fock modules ${\cal F}_c(u)$ are irreducible with respect to the Heisenberg algebra
of ${\cal E}$ generated by $\{h_r\}_{r \in {\mathbf Z}_{\neq 0}}$ with relations 
$[h_r, h_s]=\delta_{r+s,0}\frac{1}{\kappa_r}\frac{C^r-C^{-r}}{r}$. 
Let $v_c \neq 0$ be the Fock vacuum of ${\cal F}_c(u)$, we have the identification of vector spaces
$$
{\cal F}_c(u)={\mathbf C}[h_{-r}]_{r>0}v_c,~h_r v_c=0~(r>0),~C v_c=s_c v_c.
$$
The generators $e(z)$ and $\tilde{f}(z)$ are realized by vertex operators
\begin{eqnarray}
e(z) \to b_c:V_c(z;u):,~
\tilde{f}(z) \to b_c:V_c(z;u)^{-1}:,~C \to s_c,
\nonumber
\end{eqnarray}
where $b_c=-(s_c-s_c^{-1})/\kappa_1$ and
\begin{eqnarray}
V_c(z;u)=u \exp\left(\sum_{r>0}\frac{\kappa_r h_{-r}}{1-q_c^r}z^r\right)
\exp\left(\sum_{r>0}\frac{\kappa_r h_r}{1-q_c^r}q_c^{\frac{r}{2}}z^{-r}\right).\nonumber
\end{eqnarray}

\subsection{Quantum algebra ${\cal K}$}

The quantum algebra ${\cal K}$ introduced 
in Ref.\cite{Feigin-Jimbo-Mukhin-Vilkoviskiy}
is an associative algebra
with parameters $q_1, q_2, q_3$ generated by
$E_n$ $(n \in {\mathbf Z})$ and $H_r$ $(r \in {\mathbf Z}_{\neq 0})$,
and an invertible central element $C$.
We set the currents $E(z)$, $K^\pm(z)$, and $K(z)$ as follows.
\begin{eqnarray}
E(z)=\sum_{n \in {\mathbf Z}}E_nz^{-n},~
K^\pm(z)=\exp\left(\sum_{\pm r>0}H_rz^{-r}\right),~
K(z)=K^-(z)K^+(C^2z).\nonumber
\end{eqnarray}
The defining relations are as follows.
\begin{eqnarray}
&&[H_r, H_s]=-\delta_{r+s,0}
\kappa_r \frac{1+C^{2r}}{r},
\nonumber
\\
&&g(z,w)E(z)E(w)+g(w,z)E(w)E(z)\nonumber\\
&&=\frac{1}{\kappa_1}
\left(g(z,w)\delta\left(C^2\frac{z}{w}\right)K(z)+
g(w,z)\delta\left(C^2\frac{w}{z}\right)K(w)
\right),
\nonumber\\
&&
g(z,w)K^\pm(z)E(w)=\bar{g}(z,w)E(w)K^\pm(z),
\nonumber\\
&&
{\rm Sym}_{z_1, z_2, z_3}
\frac{z_2}{z_3}[E(z_1), [E(z_2), E(z_3)]]\nonumber\\
&&=
{\rm Sym}_{z_1, z_2, z_3}X(z_1, z_2, z_3)\kappa_1^{-1}
\delta\left(C^2\frac{z_1}{z_3}\right)K^-(z_1)E(z_2)K^+(z_3),
\nonumber
\end{eqnarray}
where
\begin{eqnarray}
X(z_1,z_2,z_3)&=&
\frac{(z_1+z_2)(z_3^2-z_1z_2)}{z_1z_2z_3}
G(z_2/z_3)+\frac{(z_2+z_3)(z_1^2-z_2z_3)}{z_1z_2z_3}G(z_1/z_2)
\nonumber\\
&+&\frac{(z_3+z_1)(z_2^2-z_3z_1)}{z_1z_2z_3}
\nonumber
\end{eqnarray}
and $G(w/z)$ stands for
the power series expansion of $\bar{g}(z,w)/g(z,w)$ in $w/z$.
The algebra ${\cal K}$ is a comodule over the quantum toroidal algebra ${\cal E}$.
We define the map $\Delta: {\cal K}\to {\cal E}\tilde{\otimes} {\cal K}$ as follows.
\begin{eqnarray}
&&
\Delta E(z)=e(C_2^{-1}z)\otimes K^+(z)+1\otimes E(z)+\tilde{f}(C_2z)\otimes K^-(z),
\nonumber\\
&&\Delta K^+(z)=\psi^+(C_1^{-1}C_2^{-1}z)\otimes K^+(z),
\nonumber\\
&&\Delta K^-(z)=\psi^-(C_2z)^{-1}\otimes K^-(z),
~~\Delta C=C\otimes C,\nonumber
\end{eqnarray}
where $C_1=C\otimes 1, C_2=1\otimes C$.

We introduce three families of the Fock modules ${\cal F}_{c}^B$ 
of the quantum algebra ${\cal K}$, which we call the boundary Fock modules.
We call $c$ the color.
For a complex number $s_c \in {\mathbf C}^{\times}$, let ${\cal H}_{s_c^{1/2}}$
be the Heisenberg generated by $\{H_r\}_{r \in {\mathbf Z}_{\neq 0}}$ with relations
$[H_r, H_s]=-\delta_{r+s,0}\kappa_r \frac{1+s_c^{r}}{r}$.
For $c=1,2,3$, we denote ${\cal F}_c^B$ the corresponding Fock modules of the Heisenberg algebra
${\cal H}_{s_c^{1/2}}$.
For $c=1,2,3$, the generating function $E(z)$ is realized by vertex operators
\begin{eqnarray}
E(z) \to k_c^B 
:\tilde{K}_c^-(z)\tilde{K}_c^+(s_cz):,~~C \to s_c^{1/2},
\nonumber
\end{eqnarray}
where
$k_c^B=(1+s_c)(s_d-s_b)/\kappa_1$ with $(c,d,b)=cycl(1,2,3)$,
and
\begin{eqnarray}
\tilde{K}_c^\pm (z)=\exp\left(\sum_{\pm r>0}\frac{1}{1+s_c^{-r}}H_r z^{-r}\right).
\nonumber
\end{eqnarray}

\section{Quadratic relations of ${\cal W}_{x,r}(A_{2N}^{(2)})$}
In this section, we fix a real number $r>1$ and $0<|x|<1$. 
We fix the rank $N=1,2,3,\ldots.$
Throughout this section we set
\begin{eqnarray}
q_1=x^{2r},~~q_2=x^{-2},~~q_3=x^{2(1-r)}.\nonumber
\end{eqnarray}

\subsection{Basic $W$-current}

Consider a ${\cal K}$ module defined 
as a tensor product of $N$ Fock modules ${\cal F}_2(u_i)$
of ${\cal E}$ with a boundary Fock module ${\cal F}_2^B$: 
\begin{eqnarray}
{\cal F}_2(u_1)\otimes {\cal F}_2(u_2)\otimes \cdots \otimes {\cal F}_2(u_N)\otimes {\cal F}_2^{B}.\nonumber
\end{eqnarray}
The total level is
$C=x^{-N-\frac{1}{2}}$.
The current $E(z)$ acts as a sum of vertex operators in $N+1$ bosons of the form
\begin{eqnarray}
\Delta^{(N)} E(z)=b_2 \sum_{k=1}^N \Lambda_k(z)+k_2^B \Lambda_0(z)+
b_2 \sum_{k=1}^N \Lambda_{\bar{k}}(z).\nonumber
\end{eqnarray}
Here, for $k=1,2,\ldots, N$ we set
\begin{eqnarray}
&&
\Lambda_k(z)=1\otimes \cdots \otimes V_{2}(a_k z; u_k) \otimes
\psi^+(s_2^{-1} a_{k+1}z) \otimes \cdots \otimes \psi^+(s_2^{-1}a_N z)
\otimes K^+(z),
\nonumber\\
&&\Lambda_0(z)=1\otimes \cdots \otimes 1 \otimes \tilde{K}_2(z),\nonumber\\
&&
\Lambda_{\bar{k}}(z)=1\otimes \cdots \otimes V_{2}^{-1}(a_k^{-1}z; u_k) \otimes
\psi^-(a_{k+1}^{-1}z)^{-1} \otimes \cdots \otimes \psi^-(a_N^{-1}z)^{-1}
\otimes K^-(z),
\nonumber
\end{eqnarray}
where $a_k$ are given by
$a_k=x^{N-k+\frac{1}{2}}$.
Define the dressed current 
${\mathbf \Lambda}_i(z)$ depending on 
$\mu=-x^{-2N-1}$ by
\begin{eqnarray}
{\mathbf \Lambda}_i (z)=\Lambda_i (z)
\Delta^{(N)} K_\mu^+(z)^{-1},~
K_\mu^+(z)=\prod_{s=0}^\infty K^+(\mu^{-s}z),
~i=1,\ldots,N,0,\bar{N},\ldots,\bar{1}.\nonumber
\end{eqnarray}
For $i,j=1,2,3,\ldots$ we set
\begin{eqnarray}
f_{i,j}(z)&=&\exp
\left(-\sum_{m=1}^\infty 
\frac{1}{m}[(r-1)m]_x[rm]_x(x-x^{-1})^2
{\times}\right.
\nonumber
\\&&
\times
\left.
\frac{[{\rm Min}(i,j)m]_x
\big([(N+1-{\rm Max}(i,j))m]_x-
[(N-{\rm Max}(i,j))m]_x
\big)
}{[m]_x\big(
[(N+1)m]_x-
[Nm]_x
\big)}
z^m
\right),
\nonumber
\\
d(z)&=&\frac{(1-x^{2r-1}z)(1-x^{-2r+1}z)}{(1-xz)(1-x^{-1}z)},~~~
c(x,r)=[r]_x[r-1]_x(x-x^{-1}).\nonumber
\end{eqnarray}
We define the basic $W$-current ${\mathbf T}_1(z)$ 
for ${\cal W}_{x,r}\big(A_{2N}^{(2)}\big)$ by
\begin{eqnarray}
{\mathbf T}_1(z)=
\sum_{k=1}^N {\mathbf \Lambda}_k(z)+\frac{k_2^B}{b_2} 
{\mathbf \Lambda}_0(z)+
\sum_{k=1}^N {\mathbf \Lambda}_{\bar{k}}(z).\nonumber
\end{eqnarray}
Here indices are ordered as 
$$1\prec 2\prec\cdots\prec N\prec0\prec \bar{N} \prec \cdots \prec \bar{2}\prec \bar{1}.$$
\begin{lemma}
\label{lemma1}
\cite{Feigin-Jimbo-Mukhin-Vilkoviskiy}
\footnote{
Frenkel-Reshetikhin \cite{Frenkel-Reshetikhin1}
constructed the bosonic operators 
${\mathbf \Lambda}_i^{\rm FR}(z)$ in kernel of screening operators,
that satisfy the same normal 
ordering relations as those
in Ref.\cite{Feigin-Jimbo-Mukhin-Vilkoviskiy}.}
In the case of ${\cal W}_{x,r}\big(A_{2N}^{(2)}\big)$,
the dressed currents ${\mathbf \Lambda}_i(z)$ satisfy
\begin{eqnarray}
&&f_{1,1}\left(\frac{z_2}{z_1}\right){\mathbf \Lambda}_i(z_1){\mathbf \Lambda}_j(z_2)=
d\left(\frac{x^{-1} z_2}{z_1}\right)
:{\mathbf \Lambda}_i(z_1)
{\mathbf \Lambda}_j(z_2):,~~{i \prec j, j \neq \bar{i}},
\nonumber
\\
&&f_{1,1}\left(\frac{z_2}{z_1}\right)
{\mathbf \Lambda}_j(z_1)
{\mathbf \Lambda}_i(z_2)=
d\left(\frac{x z_2}{z_1}\right)
:{\mathbf \Lambda}_j(z_1)
{\mathbf \Lambda}_i(z_2):,~~
{i \prec j, j \neq \bar{i}},
\nonumber
\\
&&f_{1,1}\left(\frac{z_2}{z_1}\right)
{\mathbf \Lambda}_0(z_1)
{\mathbf \Lambda}_0(z_2)=
d\left(\frac{z_2}{z_1}\right)
:{\mathbf \Lambda}_0(z_1){\mathbf \Lambda}_0(z_2):,
\nonumber
\\
&&f_{1,1}\left(\frac{z_2}{z_1}\right)
{\mathbf \Lambda}_i(z_1){\mathbf \Lambda}_i(z_2)=
:{\mathbf \Lambda}_i(z_1){\mathbf \Lambda}_i(z_2):,~~{i \neq 0},
\nonumber
\\
&&f_{1,1}\left(\frac{z_2}{z_1}\right)
{\mathbf \Lambda}_k(z_1){\mathbf \Lambda}_{\bar{k}}(z_2)=
d\left(\frac{x^{-1} z_2}{z_1}\right)
d\left(\frac{x^{-2N-2+2k}z_2}{z_1}\right)
:{\mathbf \Lambda}_k(z_1){\mathbf \Lambda}_{\bar{k}}(z_2):,
~{1\leq k \leq N},
\nonumber
\\
&&f_{1,1}\left(\frac{z_2}{z_1}\right){\mathbf \Lambda}_{\bar{k}}(z_1){\mathbf \Lambda}_{k}(z_2)=
d\left(\frac{x z_2}{z_1}\right)
d\left(\frac{x^{2N+2-2k}z_2}{z_1}\right)
:{\mathbf \Lambda}_{\bar{k}}(z_1)
{\mathbf \Lambda}_{k}(z_2):,~{1\leq k \leq N}.
\nonumber
\end{eqnarray}
\end{lemma}
\subsection{Quadratic relations}
In this section, we introduce the higher $W$-currents
${\mathbf T}_j(z)$ and obtain quadratic relations of them.
We define the higher $W$-currents ${\mathbf T}_i(z)$ 
$(i \in {\mathbf N})$ for ${\cal W}_{x,r}\big(A_{2N}^{(2)}\big)$
by fusion relation
\begin{eqnarray}
&&
\lim_{z_1 \to x^{\pm (i+j)}z_2}
\left(1-\frac{x^{\pm (i+j)}z_2}{z_1}\right)f_{i,j}\left(\frac{z_2}{z_1}\right){\mathbf T}_i(z_1){\mathbf T}_j(z_2)\nonumber
\\
&=&
\mp c(x,r)\prod_{l=1}^{{\rm Min}(i,j)-1}d(x^{2l+1})
~{\mathbf T}_{i+j}(x^{\pm i}z_2),~
{i, j \geq 1},
\nonumber
\end{eqnarray}
and ${\mathbf T}_0(z)=1$.
We obtain ${\mathbf T}_i(z)=0$ for $ i \geq 2N+2$.
\begin{proposition}
\cite{Kojima3}
In the case of
${\cal W}_{x,r}\big(A_{2N}^{(2)}\big)$, 
the $W$-currents ${\mathbf T}_i(z)$ satisfy the duality
\begin{eqnarray}
{\mathbf T}_{2N+1-i}(z)=\frac{[r-\frac{1}{2}]_x}{[\frac{1}{2}]_x}\prod_{k=1}^{N-i}
d(x^{2k}){\mathbf T}_i(z),~~0\leq i \leq N.
\label{duality}
\end{eqnarray}
\end{proposition}
\begin{theorem}
\cite{Kojima3}
In the case of ${\cal W}_{x,r}\big(A_{2N}^{(2)}\big)$,
the $W$-currents ${\mathbf T}_i(z)$
satisfy the set of quadratic relations
\begin{eqnarray}
&&
f_{i,j}\left(\frac{z_2}{z_1}\right){\mathbf T}_i(z_1)
{\mathbf T}_j(z_2)-
f_{j,i}\left(\frac{z_1}{z_2}\right)
{\mathbf T}_j(z_2){\mathbf T}_i(z_1)
\nonumber\\
&=& c(x,r)\sum_{k=1}^i \prod_{l=1}^{k-1}d(x^{2l+1})\times
\nonumber\\
&&\times
\left(
\delta\left(\frac{x^{-j+i-2k}z_2}{z_1}\right)
f_{i-k, j+k}(x^{j-i}){\mathbf T}_{i-k}(x^{k}z_1)
{\mathbf T}_{j+k}(x^{-k}z_2)
\right.\nonumber
\\
&&-\left.\delta\left(\frac{x^{j-i+2k}z_2}{z_1}\right)
f_{i-k, j+k}(x^{-j+i}){\mathbf T}_{i-k}(x^{-k}z_1)
{\mathbf T}_{j+k}(x^kz_2)
\right)\nonumber\\
&&+
c(x,r)\prod_{l=1}^{i-1}
d(x^{2l+1})
\prod_{l=N+1-j}^{N+i-j}d(x^{2N})
\nonumber\\
&&
{\times}
\left(
\delta\left(\frac{x^{-2N+j-i-1}z_2}{z_1}\right)
{\mathbf T}_{j-i}(x^{-i}z_2)
-
\delta\left(\frac{x^{2N-j+i+1}z_2}{z_1}\right)
{\mathbf T}_{j-i}(x^{i}z_2)
\right),\nonumber\\
&&~~~1\leq i \leq j \leq N.
\label{quadratic1}
\end{eqnarray}
\end{theorem}
\begin{definition}
Let $W$ be the free complex associative algebra generated by elements 
$\overline{\mathbf T}_i[m], m\in {\mathbf Z}, 0\leq i \leq 2N+1$, 
$I_K$ the left
ideal generated by elements $\overline{\mathbf T}_i[m], m \geq K \in {\mathbf N}, 
0\leq i \leq 2N+1$, and
\begin{eqnarray}
\widehat{W}={\displaystyle \lim_{\leftarrow}W/I_K}.
\nonumber
\end{eqnarray}
The deformed $W$-algebra 
${\cal W}_{x,r}\bigl(A_{2N}^{(2)}\bigr)$
is the quotient of $\widehat{W}$ by the two-sided ideal generated by
the coefficients of the generating series which are the differences
of the right hand sides and of the left hand sides
of the relations (\ref{duality}) and (\ref{quadratic1}),
where the generating series ${\mathbf T}_i(z)$
are replaced with
${
\overline{\mathbf T}_i(z)=\sum_{m \in {\mathbf Z}}
\overline{\mathbf T}_i[m]z^{-m}}, 0\leq i \leq 2N+1$.
\end{definition}

\section{Quadratic relations of $q$-$Y_{L_1, L_2, L_3}$}

We fix natural numbers $L_1, L_2, L_3$
such that 
$L_1+L_2+L_3\geq 1$, $L_1, L_2, L_3
\in {\mathbf N}$.
We fix real numbers $\lambda_1, \lambda_2, \lambda_3$ such that
$\lambda_1+\lambda_2+\lambda_3=0$.
We fix $0<|x|<1$.
Throughout this section we set
\begin{eqnarray}
q_1=x^{2\lambda_1},~q_2=x^{2\lambda_2},
~q_3=x^{2\lambda_3},
~L=L_1+L_2+L_3.\nonumber
\end{eqnarray}
\subsection{Basic $W$-current}
Consider a tensor product of $L$ Fock modules of ${\cal E}$:
\begin{eqnarray}
{\cal F}_{c_1}(u_1)\otimes \cdots \otimes {\cal F}_{c_L}(u_L).\nonumber
\end{eqnarray}
Here we choose colors $c_1, c_2,\ldots, c_{L}\in \{1,2,3\}$
such that 
$$L_c=n\big(I(c)\big),
~~I(c)=\{1\leq i \leq L|c_i=c\}~~(c=1,2,3).$$
The total level of $L$ Fock modules is
$C=\prod_{j=1}^L s_{c_j}$.
The current $e(z)$ acts as a sum of vertex operators in $L$ 
bosons of the form
\begin{eqnarray}
\Delta^{(L-1)} e(z)=\sum_{i=1}^L b_{c_i}\Lambda_i(z).
\nonumber
\end{eqnarray}
Here, for $i=1,2,\ldots, L$ we set
\begin{eqnarray}
&&
\Lambda_i(z)=1\otimes \cdots \otimes 1 \otimes
V_{c_i}(a_i z; u_i) \otimes
\psi^+(s_{c_{i+1}}^{-1} a_{i+1}z) \otimes \cdots \otimes 
\psi^+(s_{c_L}^{-1}a_{L} z),
\nonumber
\end{eqnarray}
where $a_i$ are given by
$a_i=\prod_{j=i+1}^L s_{c_j}^{-1}$.
Define the dressed currents 
${\mathbf \Lambda}_i(z)$ depending on
free parameter $\mu$ by
\begin{eqnarray}
{\mathbf \Lambda}_i (z)=\Lambda_i (z)
\Delta^{(L-1)} K_\mu^+(z)^{-1},~
K_\mu^+(z)=\prod_{s=0}^\infty K^+(\mu^{-s}z),
~\mu=x^{2\alpha}.\nonumber
\end{eqnarray}
For $i,j=1,2,3,\ldots$, we set
\begin{eqnarray}
g_{i,j}(z)&=&\exp
\left(\sum_{m=1}^\infty 
\frac{1}{m}
\frac{
[\lambda_1 m]_x[\lambda_3 m]_x
[{\rm Min}(i,j)\lambda_2 m]_x
[(\alpha-{\rm Max}(i,j)\lambda_2)m]_x
}{[\lambda_2 m]_x[\alpha m]_x}(x-x^{-1})^2
z^m
\right),
\nonumber
\\
d_\lambda(z)
&=&
\frac{(1-x^{\lambda_1-\lambda_3}z)
(1-x^{\lambda_3-\lambda_1}z)}{
(1-x^{\lambda_2}z)(1-x^{-\lambda_2}z)},~~~
c_\lambda=\frac{[\lambda_1]_x[\lambda_3]_x}{[\lambda_2]_x}(x-x^{-1}),
\nonumber
\\
\gamma_c(z)
&=&
\frac{(1-x^{2\lambda_c}z)(1-x^{-2\lambda_c}z)}{
(1-x^{2\lambda_2}z)(1-x^{-2\lambda_2}z)},~~~c=1,2,3.
\nonumber
\end{eqnarray}
We define the basic $W$-current ${\mathbf T}_1(z)$ 
for $q$-$Y_{L_1, L_2, L_3}$ by
\begin{eqnarray}
{\mathbf T}_1(z)=
\sum_{i=1}^L b_{c_i}{\mathbf \Lambda}_i(z).
\nonumber
\end{eqnarray}
\begin{lemma}
\cite{Feigin-Jimbo-Mukhin-Vilkoviskiy}
In the case of $q$-$Y_{L_1, L_2, L_3}$,
the dressed currents ${\mathbf \Lambda}_i(z)$ satisfy
\begin{eqnarray}
&&
g_{1,1}\left(\frac{z_2}{z_1}\right)
{\mathbf \Lambda}_i(z_1){\mathbf \Lambda}_j(z_2)=
d_\lambda\left(\frac{x^{\lambda_2} z_2}{z_1}\right)
:{\mathbf \Lambda}_i(z_1)
{\mathbf \Lambda}_j(z_2):,~~{1\leq i<j \leq L},
\nonumber
\\
&&
g_{1,1}\left(\frac{z_2}{z_1}\right)
{\mathbf \Lambda}_j(z_1)
{\mathbf \Lambda}_i(z_2)=
d_\lambda\left(\frac{x^{-\lambda_2} z_2}{z_1}\right)
:{\mathbf \Lambda}_j(z_1)
{\mathbf \Lambda}_i(z_2):,~~
{1\leq i<j \leq L},
\nonumber
\\
&&
g_{1,1}\left(\frac{z_2}{z_1}\right)
{\mathbf \Lambda}_i(z_1)
{\mathbf \Lambda}_i(z_2)=
\gamma_{c_i}\left(\frac{z_2}{z_1}\right)
:{\mathbf \Lambda}_i(z_1){\mathbf \Lambda}_i(z_2):,~~1\leq i \leq L.
\nonumber
\end{eqnarray}
\end{lemma}
\subsection{Quadratic relations}
\label{sec:3}
In this section, we introduce the higher $W$-currents
${\mathbf T}_i(z)$ and obtain quadratic relations of them.
We define the higher $W$-currents ${\mathbf T}_i(z)$, $(i \in {\mathbf N})$
for $q$-$Y_{L_1,L_2,L_3}$
by fusion relation
\begin{eqnarray}
&&
\lim_{z_1 \to x^{\pm (i+j)\lambda_2}z_2}
\left(1-\frac{x^{\pm (i+j)\lambda_2}z_2}{z_1}\right)
g_{i,j}\left(\frac{z_2}{z_1}\right)
{\mathbf T}_i(z_1){\mathbf T}_j(z_2)\nonumber
\\
&=&
\mp c_\lambda \prod_{l=1}^{{\rm Min}(i,j)-1}
d_\lambda(x^{(2l+1)\lambda_2})
~{\mathbf T}_{i+j}(x^{\pm i \lambda_2}z_2),~
{i, j \geq 1},
\nonumber
\end{eqnarray}
and ${\mathbf T}_0(z)=1$.
\begin{proposition}
In the case of $q$-$Y_{L_1,L_2,L_3}$, 
the $W$-currents ${\mathbf T}_i(z)$ don't vanish.
\begin{eqnarray}
{\mathbf T}_i(z)\neq 0,~~~i\in {\mathbf N}.\nonumber
\end{eqnarray}
\end{proposition}

Upon the specialization
$(L_1,L_2,L_3)=(0, N+1, 0)$,
the $W$-currents ${\mathbf T}_i(z)$ of $q$-$Y_{L_1,L_2,L_3}$ coincide
with those of
${\cal W}_{x,r}\big(A_N^{(1)}\big)$ in Refs.\cite{
Shiraishi-Kubo-Awata-Odake, 
Awata-Kubo-Odake-Shiraishi}.
Upon the specialization
$(L_1,L_2,L_3)=(0, N+1, M+1)$,
the $W$-currents ${\mathbf T}_i(z)$ of $q$-$Y_{L_1,L_2,L_3}$
coincide 
with those of
${\cal W}_{x,r}\big(A(M, N)^{(1)}\big)$ in Refs.\cite{
Kojima1, Kojima2}.
\begin{eqnarray}
&&
{\mathbf T}_i(z)=0~~(i \geq N+1)~~~{\rm for}
~~{\cal W}_{x,r}\big(A_N^{(1)}\big):~\textrm{non-super},\nonumber\\
&&
{\mathbf T}_i(z)\neq 0~~(i \in {\mathbf N})~~~{\rm for}
~~{\cal W}_{x,r}\big(A(M,N)^{(1)}\big):~{\rm super}.\nonumber
\end{eqnarray}
\begin{theorem}
\label{theorem2}
\footnote{
Harada, Matsuo, Noshita and Watanabe
conjectured similar quadratic relations
in Ref.\cite{Harada-Matsuo-Noshita-Watanabe}.
The quadratic relations (\ref{quadratic2}) of
$q$-$Y_{L_1, L_2, L_3}$ in Theorem \ref{theorem2} 
can be proved similarly to those of
${\cal W}_{x,r}\big(A(M,N)^{(1)}\big)$ 
in Refs.\cite{Kojima1, Kojima2}.}
In the case of $q$-$Y_{L_1, L_2, L_3}$,
the $W$-currents ${\mathbf T}_i(z)$
satisfy the set of quadratic relations
\begin{eqnarray}
&&
g_{i,j}\left(\frac{z_2}{z_1}\right){\mathbf T}_i(z_1)
{\mathbf T}_j(z_2)-
g_{j,i}\left(\frac{z_1}{z_2}\right)
{\mathbf T}_j(z_2){\mathbf T}_i(z_1)\nonumber
\\
&=& c_\lambda
\sum_{k=1}^i \prod_{l=1}^{k-1}d_\lambda(x^{(2l+1)\lambda_2})\times
\nonumber\\
&&
\times
\left(
\delta\left(\frac{x^{(j-i+2k)\lambda_2}z_2}{z_1}\right)
g_{i-k, j+k}(x^{(i-j)\lambda_2}){\mathbf T}_{i-k}(x^{-k\lambda_2}z_1)
{\mathbf T}_{j+k}(x^{k\lambda_2}z_2)
\right.
\nonumber
\\
&&-\left.\delta\left(\frac{x^{(-j+i-2k)\lambda_2}z_2}{z_1}\right)
g_{i-k, j+k}(x^{(j-i)\lambda_2}){\mathbf T}_{i-k}(x^{k\lambda_2}z_1)
{\mathbf T}_{j+k}(x^{-k\lambda_2}z_2)
\right),\nonumber
\\
&&~~~1\leq i \leq j.
\label{quadratic2}
\end{eqnarray}
\end{theorem}

Upon the specialization
$(L_1,L_2,L_3)=(0, N+1, M+1)$ and
$(\lambda_1,\lambda_2,\lambda_3)=(r,-1,1-r)$,
the set of quadratic relations (\ref{quadratic2}) 
of $q$-$Y_{L_1, L_2, L_3}$ 
coincides with those of
${\cal W}_{x,r}\big(A(M, N)^{(1)}\big)$ in Refs.\cite{Kojima1, Kojima2}.
~\\
\begin{definition}
Let $W$ be the free complex associative algebra generated by elements 
$\overline{\mathbf T}_i[m], m\in {\mathbf Z}, i\in {\mathbf N}$, 
$I_K$ the left
ideal generated by elements $\overline{\mathbf T}_i[m], m \geq K \in {\mathbf N}, 
i \in {\mathbf N}$, and
\begin{eqnarray}
\widehat{W}={\displaystyle \lim_{\leftarrow}W/I_K}.
\nonumber
\end{eqnarray}
The $q$-$Y_{L_1, L_2, L_3}$
is the quotient of $\widehat{W}$ by the two-sided ideal generated by
the coefficients of the generating series which are the differences
of the right hand sides and of the left hand sides
of the relations (\ref{quadratic2}),
where the generating series ${\mathbf T}_i(z)$
are replaced with
${
\overline{\mathbf T}_i(z)=\sum_{m \in {\mathbf Z}}
\overline{\mathbf T}_i[m]z^{-m}}, i\in {\mathbf N}$.
\end{definition}

\begin{acknowledgement}
The author is thankful for the kind hospitality by the organizing committee of the 15-th International Workshop ''Lie theory and its application in physics'' (LT15).
The author would like to thank Professor Michio Jimbo for giving advice.
This work is supported by 
the Grant-in-Aid for Scientific Research {\bf C}(19K0350900)
from Japan Society for Promotion of Science.
\end{acknowledgement}
\end{document}